\documentclass[twoside,leqno,twocolumn]{article}

\usepackage[letterpaper]{geometry}

\usepackage{ltexpprt}
\usepackage{hyperref}
\usepackage{algpseudocode}
\usepackage{algorithmicx}
\usepackage{algorithm}
\usepackage{setspace}
\usepackage{xcolor}
\usepackage{subcaption}
\usepackage{graphicx}
\usepackage{wrapfig}
\usepackage{amsmath}
\usepackage{amssymb}

\usepackage{siunitx} %

\algrenewcommand\textproc{}

\usepackage[capitalise,noabbrev]{cleveref}
\crefformat{equation}{(#2#1#3)}

\clearpage

\newcommand{\G}{\mathcal{G}}
\newcommand{\F}{\mathcal{F}}
\newcommand{\BigO}{\mathcal{O}}

\newcommand{\Dgione}{\Delta g_{i-1}}

\newcommand{\Dfione}{\Delta f_{i-1}}

\clearpage

\newcommand{\code}[1]{\texttt{#1}}
\renewcommand{\vec}[1]{\boldsymbol{\mathrm{#1}}}

\begin{document}

\title{\Large Performance of Low Synchronization Orthogonalization Methods in Anderson Accelerated Fixed Point Solvers
}
\author{Shelby Lockhart\thanks{University of Illinois at Urbana-Champaign}
\and David J. Gardner\thanks{Lawrence Livermore National Laboratory}
\and Carol S. Woodward\footnotemark[2]
\and Stephen Thomas\thanks{National Renewable Energy Laboratory}
\and Luke N. Olson\footnotemark[1]
}

\date{}

\maketitle

\maketitle

\begin{abstract} {\small\baselineskip=9pt
Anderson Acceleration (AA) is a method to accelerate 
the convergence of fixed point
iterations for nonlinear, algebraic systems of equations.  Due to the requirement of solving a least squares problem at each iteration and a 
reliance on modified Gram-Schmidt for updating the iteration space, AA requires extra costly synchronization steps for global reductions. Moreover, the number of reductions in each iteration depends on the size of the iteration space. 
In this work, we introduce three low synchronization orthogonalization algorithms into AA within SUNDIALS that reduce the total number of global reductions per
iteration to a constant of 2 or 3, independent of the size of the iteration space.
A performance study demonstrates the reduced time required by the new
algorithms at large processor counts with CPUs
and demonstrates the predicted performance on multi-GPU architectures.
Most importantly, we provide convergence and timing data for multiple numerical
experiments to demonstrate reliability of the algorithms within AA and improved performance at parallel strong-scaling limits.
}\end{abstract}

\section{Introduction}
Anderson acceleration (AA) is a method employed to
accelerate the convergence of the fixed point (FP) method for solving systems of nonlinear equations \cite{anderson1965iterative,WALKER_NI,kelley2018numerical}. 
The formulation of this method leading to the most efficient implementation requires 
the solution of an unconstrained minimization problem, generally done through solving a  
least squares problem (LSP) via $QR$ factorization \cite{WALKER_NI,Fang-Saad09}. Because
the LSP is highly sensitive, it is imperative to employ a stable algorithm,
such as modified Gram-Schmidt, to update the factorization at each iteration. The
application of modified Gram-Schmidt generates the main bottleneck of the algorithm
when performed in a distributed memory parallel environment, due to the associated
high communication costs of performing multiple dot products in sequence
\cite{fraysse1998}.

A solution to the costly $QR$ updates within AA is to apply an
 algorithm that requires fewer synchronization points per iteration, i.e.
performs multiple dot products simultaneously or computes an equivalent
$QR$ update at a lower parallel cost. 
Recent research by \'{Swirydowicz} et al.\ \cite{LOWSYNC_GMRES} and Bielich et al.\ \cite{DCGS2Arnoldi}
introduced low synchronization algorithms for
MGS, CGS-2 and GMRES \cite{Saad86}, 
and because AA is equivalent to GMRES when solving systems of
linear equations \cite{WALKER_NI}, it is natural to consider leveraging these advances in AA.

We have implemented and tested three low synchronization algorithms within the SUNDIALS KINSOL \cite{KINSOL} implementation of AA. 
These methods reduce the number of global reductions to a constant number per iteration independent of the size of the AA iteration space.
While these methods can reduce the number of synchronizations to 1 within  Krylov methods, a single reduction is not possible within AA due to the requirement for normalization before the LSP solve.  Nevertheless, these
new methods are able to reduce required global communications to 2 or 3 per iteration.

In addition to demonstrating how to extend these ideas to AA, this paper also
includes a problem independent performance study that illustrates the improved strong scalability of the new kernels, particularly for large iteration spaces and processor counts.  The relative cost of the new methods is dependent upon on particular problem. Thus, we present convergence and timing data for multiple examples to demonstrate the reliability of the algorithms for AA and improved performance at scale.

The rest of the paper is organized as follows.  In \cref{sec:background}, we describe the AA algorithm and relevant components in the $QR$ factorization process.  
In \cref{sec:qradd}, we present and compare the low synchronization $QR$ factorization variants. In \cref{sec:perf_study,sec:experiments} we give performance results for the different methods employed within the AA algorithm itself and from several test problems respectively. Finally, \cref{sec:conclusions} provides conclusions based on our findings.

\section{Background}
\label{sec:background}

\subsection{Anderson Acceleration}
\label{ss:anderson}

AA (\cref{alg:anderson}) utilizes a linear combination of $m$ prior function evaluations to accelerate the convergence of the FP iteration $x_{i+1} = G(x_i)$. The weights $\gamma_k^{(i)}$  minimize the FP residual norm in the linear case and are computed by solving a LSP (\cref{LSP} in \cref{alg:anderson}) via a $QR$ factorization of $\F_i$ where $Q$ is a matrix with linearly
independent, normalized columns, and $R$ is an upper triangular matrix.
\begin{center}
\vspace{-2em}
\resizebox{.5\textwidth}{!}{%
\begin{minipage}{0.575\textwidth}
\begin{algorithm}[H]
    \caption{Anderson Acceleration (AA)}\label{alg:anderson}
    \hspace*{\algorithmicindent} \textbf{Input:} $x_0$, $m \geq 1$, tol, and maxIters. \\
    \hspace*{\algorithmicindent} \textbf{Output:} $x^*$
    \begin{algorithmic}[1]
        \State $x_1 = G(x_0)$ and $f_0 = G(x_0) - x_0$
        \For{$i = 1,2,\ldots,\text{maxIters}$}
            \State Set $m_i = \text{min}\{m,i\}$ and $f_i = G(x_i) - x_i$
            \State Set $\G_i = [\Delta g_{i-m_i},\ldots, \Dgione]$ with $\Delta g_k = G(x_{k+1}) - G(x_{k})$
            \State Set $\F_i = [\Delta f_{i-m_i},\ldots, \Dfione]$ with $\Delta f_k = f_{k+1} - f_{k}$
            \State \label{LSP}Determine $\gamma^{(i)} = [\gamma_0^{(i)},\ldots, \gamma_{m_{i-1}}^{(i)}]^T$ that solves $$\text{min}_{\gamma^{(i)}} \| f_i - \F_i \gamma^{(i)} \|_2.$$
            \State Set $x_{i+1} = G(x_i) - \G_i \gamma^{(i)}$
            \If{$\| x_{i+1} - x_{i} \| < \text{tol}$}
                \State \Return $x_{i+1}$ as $x^*$
            \EndIf
        \EndFor
        \end{algorithmic}
\end{algorithm}
\end{minipage}
}
\end{center}
Solving the LSP (\cref{alg:qr_lsp}) requires updating the $QR$ decomposition to incorporate the vector $\Dfione$. This update is comprised of two subroutines: \code{QRAdd} appends a vector to the $QR$ factorization while performing the necessary orthogonalization and \code{QRDelete} removes the oldest vector from the factorization.
\begin{center}
\vspace{-2em}
\resizebox{.5\textwidth}{!}{%
\begin{minipage}{0.575\textwidth}
\begin{algorithm}[H]
    \caption{Least Squares Problem (LSP) Solve}\label{alg:qr_lsp}
    \hspace*{\algorithmicindent} \textbf{Input:} $i$, $m$, $m_i$, $f_i$, $\Dfione$, $\gamma^{(i)}$, $Q$, and $R$. \\
    \hspace*{\algorithmicindent} \textbf{Output:} $Q$, $R$, and $\gamma^{(i)}$
    \begin{algorithmic}[1]
        \If{$i = 1$}
            \State Set $Q_{:,0} = \Dfione / \|\Dfione \|_2$ and $R_{0,0} = \| \Dfione \|_2$
        \Else
            \If{$i > m$}
                \State \Call{\code{QRDelete}}{$Q$, $R$, $m_i$} \label{anderson_qrdelete}
            \EndIf
            \State \Call{\code{QRAdd}}{$Q$, $R$, $\Dfione$, $m_i$} \label{anderson_qradd}
        \EndIf
        \State Solve $R \gamma^{(i)} = Q^T * f_i$ \label{qr_solve}   
        \end{algorithmic}
\end{algorithm}
\end{minipage}
}
\end{center}
A key factor in the performance and convergence of AA is the number of residual history vectors or depth, $m$. The depth determines the maximum size of the $QR$ factorization and thus the number of vectors that must be orthogonalized each iteration in \code{QRAdd}. Note in the start-up phase ($i < m$), the number of vectors increases by one each iteration until $m$ vectors are retained. Once the history is full, the subsequent iterations require orthogonalization against $m - 1$ vectors. As discussed below, the number of synchronizations required by \code{QRAdd} depends on $m$ and the method for updating the $QR$ factorization. Generally, \code{QRDelete} does 
not require any communication, as initially noted by
Loffeld and Woodward in \cite{loffeldwoodwardconsiderations}. In practice the depth is often kept small ($m < 5$) however, there are potential convergence benefits to be gained if it were economical, performance-wise, to run with larger $m$. 

\subsection{QR Factorization with Gram-Schmidt}
\label{ss:qrfac}

There exist several methods for computing the $QR$ factorization of
a matrix. In this paper, we consider methods derived from 
the Gram-Schmidt procedure. Specifically, we focus on the process
of updating a given $QR$ factorization to include a vector $q$.
Here, column-oriented algorithms are presented as opposed to
those that proceed row-wise. These are better suited to a
low-synchronization formulation and promote data-locality for
cache-memory access both for CPUs and GPUs.

The classical Gram-Schmidt (CGS) process updates a given $QR$ factorization
with the vector $q$ by applying the projection $P = I - QQ^T$ to $q$ then normalizing
$q$. While this process is attractive in a parallel computing environment,
due to the fact that $P$ can be applied with a single
reduction, this algorithm is highly unstable and
exhibits a loss of orthogonality dependent on $\BigO (\varepsilon) \kappa^2(A)$
\cite{giraud20051, giraud20052} where $\varepsilon$ is machine precision and $\kappa^2(A)$ is the square of
the condition number of $A$. For this reason alone, 
it is unsuitable as a $QR$ updating
scheme in AA, where the LSP solver convergence may be highly sensitive
to the loss of orthogonality \cite{paige2018-Book-}. 

A common solution for the instability of CGS is to apply modified Gram-Schmidt (MGS)
instead. This process applies the rank--1 elementary projection matrices 
$I - Q_{:,j}Q_{:,j}^T$, $j=0,\ldots,n-1$, one for each of the $n$
linearly independent columns of $Q$, to the vector $q$ then normalizes $q$. 
Performing the Gram-Schmidt process
in this way effectively reduces the loss of orthogonality to
$\BigO(\varepsilon) \kappa (A)$ \cite{bjorck1967}. However, applying the projections sequentially requires $n$ dot products, 
one for each of the columns, resulting in
communication bottlenecks when performed in a parallel 
distributed environment (discussed more thoroughly in \cref{ss:qradd_mgs}).

New variants of CGS and MGS have been introduced to mitigate instability and
parallel performance issues caused by these methods. In most cases, they
introduce a correction matrix $T$, leading to a projection operator of the form $P = I - Q T Q^T$. 
In the case of CGS, CGS-2 (classical Gram-Schmidt with reorthogonalization)
corrects the projection by reorthogonalizing the vectors of $Q$ and thereby
reduces the loss of orthogonality to $\BigO (\varepsilon)$ \cite{giraud20052, parlett1998}.
The form of the correction matrix for this algorithm was derived
in \cite{LOWSYNC_GMRES} and is discussed further in \cref{ss:qradd_cgs2,ss:qradd_dcgs2}.
The inverse compact $WY$ MGS representation was recently introduced,
requiring a triangular solve instead of the recursive construction of the
correction matrix $T$ from a block-triangular inverse
(presented in \cref{ss:qradd_icwy}).
The required number of dot products per iteration is effectively reduced 
to one for GMRES and s--step Krylov solvers
\cite{LOWSYNC_GMRES, lowsync_sstep}. The inverse compact $WY$
algorithm maintains  $\BigO (\varepsilon) \kappa (A)$ loss of orthogonality.

\section{QR Update Methods in AA}
\label{sec:qradd}
In this section, we discuss the \code{QRAdd} kernel, 
\cref{anderson_qradd} of \cref{alg:qr_lsp},
employed to update the AA iteration space with the
$\Dfione$ vector. 
We present the baseline \code{QRAdd\_MGS}
(modified Gram Schmidt) algorithm alongside
three low synchronization variants of the kernel implemented
within AA in SUNDIALS.
For each of the \code{QRAdd} kernels, we discuss
the form of the projectors applied to orthogonalize the given set of
vectors, $\F_i$, as well as, their predicted parallel performance.

Throughout this section and the remainder of the paper, 
matrix entries are referenced via subscripts with 0-based indexing. 
For example,
$M_{0,0}$ refers to the first row, first column of a matrix $M$,
and slices of a matrix are inclusive, i.e. 
$M_{:,0:k-1}$ refers to all rows of the matrix and columns $0,\ldots,k-1$.

\subsection{Modified Gram Schmidt}\label{ss:qradd_mgs}
The standard approach for updating the $QR$ factorization within AA
is to apply MGS, outlined in \cref{alg:qr_mgs}.
\begin{center}
\vspace{-2em}
\resizebox{.5\textwidth}{!}{%
\begin{minipage}{0.575\textwidth}
\begin{algorithm}[H]
    \setstretch{1.15}
    \caption{\code{QRAdd\_MGS}}\label{alg:qr_mgs}
    \hspace*{\algorithmicindent} \textbf{Input:} $Q$, $R$, $\Dfione$ and $m_i$ \\
    \hspace*{\algorithmicindent} \textbf{Output:} $Q$, $R$
    \begin{algorithmic}[1]
        \For{$j=0,\ldots,m_i-2$}
                \State $R_{j,m_i-1} \leftarrow Q_{:,j}^T * \Dfione$
                \Comment{Sync}
            \State $\Dfione \leftarrow 
                    \Dfione - R_{j,m_i-1} * Q_{:,j}$\label{mgs_linearsum}
        \EndFor
        \State $R_{m_i-1,m_i-1} \leftarrow \| \Dfione \|_2$
        \Comment{Sync}
        \State $Q_{:,m_i-1} \leftarrow \Dfione / R_{m_i-1,m_i-1}$
    \end{algorithmic}
\end{algorithm}
\end{minipage}
}
\end{center}
In each AA iteration, applying MGS requires $m_i$ dot products; $m_i - 1$ for the orthogonalization against all previous vectors in the AA
iteration space, and one for the normalization at the end of the algorithm. This algorithm results
in $m_i$ synchronizations across all processes, as the reductions form
a chain of dependencies and can not be performed in tandem. 
The high costs of these synchronizations are exacerbated by the lack of
computational workload between each reduction, namely the entire
algorithm only requires $\BigO (m_i n)$ flops.

\subsection{ICWY Modified Gram Schmidt}\label{ss:qradd_icwy}

Low-synch one-reduce Gram-Schmidt algorithms are based upon
two key ideas. First, the compact $WY$ representation relies on
a triangular correction matrix $T$, which contains a 
strictly lower triangular matrix $L$. 
One row or block of rows of $L$ is computed at a time in a 
single global reduction.
Each row 
$$L_{m_i-2,0:m_i-2} = (Q_{:,0:m_i-2}^T Q_{m_i-2})^T$$
is obtained within the current step, then the
normalization step is lagged and merged into a single reduction.
The associated orthogonal projector is based on %
R\"uhe \cite{ruhe1983} and presented in \'{S}wiryzdowicz  et al.\
\cite{LOWSYNC_GMRES}.
\begin{align*}
P & = I - Q_{0:m_i-2} T_{0:m_i-2,0:m_i-2} Q_{:,0:m_i-2}^T, \\
T_{0:m_i-2,0:m_i-2} & = (\: I + L_{m_i-1} \:)^{-1} \\
                    & \approx (Q_{:, 0:m_i-2}^T Q_{:, 0:m_i-2})^{-1}.
\end{align*}
The implied triangular solve requires an additional $(m_i-1)^2$ flops at
iteration $m_i-1$ and thus leads to a slightly higher operation count 
compared to the original MGS algorithm. 
The operation $Q_{:,0:m_i-2}^T Q_{:,m_i-2}$ increases 
ICWY-MGS complexity by $m_in^2$ for an overall complexity of $\BigO (m_i n^2)$,
and reduces synchronizations from $m_i-1$ at iteration 
$m_i$ to 1. 
Only one global reduction is required per iteration, and the 
amount of inter-process communication does not depend upon the number  
of rank--1 projections $I - Q_{:,j}\:Q_{:,j}^T$ applied at each iteration.

When implementing ICWY-MGS in the context of AA, lagging the
normalization until the subsequent iteration is not an option,
as the factorization is applied immediately after updating in the 
LSP solve on \cref{LSP} of \cref{alg:qr_lsp}. The resulting
$QR$ update algorithm within AA is detailed in \cref{alg:qr_icwy_mgs}.

\begin{center}
\vspace{-2em}
\resizebox{.49\textwidth}{!}{%
\begin{minipage}{0.575\textwidth}
\begin{algorithm}[H]
    \setstretch{1.15}
    \caption{\code{QRAdd\_ICWY}}\label{alg:qr_icwy_mgs}
    \hspace*{\algorithmicindent} \textbf{Input:} $Q$, $R$, $T$, $\Dfione$ and $m_i$ \\
    \hspace*{\algorithmicindent} \textbf{Output:} $Q$, $R$, $T$
    \begin{algorithmic}[1]
        \State $T_{m_i-2,0:m_i-2} \leftarrow Q_{:,0:m_i-2}^T * Q_{:,m_i-2}$
        \Comment{Delayed Sync}\label{icwy_delayedsync}
        \State $R_{0:m_i-2,m_i-1} \leftarrow Q_{:,0:m_i-2}^T * \Dfione$
        \Comment{Sync}\label{icwy_sync}
        \State $T_{m_i-2,m_i-2} \leftarrow 1$
        \State $R_{0:m_i-2, m_i-1} \leftarrow T_{0:m_i-2,0:m_i-2}^{-1}R_{0:m_i-2,m_i-1}$
        \State $\Dfione \leftarrow \Dfione - Q_{:,0:m_i-2} * R_{0:m_i-2,m_i-1}$
            \label{icwy_linearsum}
        \State $R_{m_i-1,m_i-1} \leftarrow \| \Dfione \|_2$
        \Comment{Sync}
        \State $Q_{:,m_i-1} \leftarrow \Dfione / R_{m_i-1,m_i-1}$
    \end{algorithmic}
\end{algorithm}
\end{minipage}
}
\end{center}
Here, we present a two reduction variant of the ICWY-MGS algorithm. 
The formation of the correction operator and the matrix $R$ can
be merged on \cref{icwy_sync}. With the inclusion of the normalization
at the end of the algorithm, this results in two synchronizations per
iteration until the AA iteration space is filled.
Once AA reaches \cref{anderson_qrdelete} in \cref{alg:qr_lsp},
the oldest vector in the factorization is deleted. 
While this can be performed with Givens rotations in cases
where only $Q$ and $R$ are stored,
the explicit storage and application of $T$ requires that this correction
matrix be updated by introducing a single reduction.
Overall, this process results in two global reduction steps per iteration
until $m$ is reached, after which, there are three global reductions per
iteration.

\subsection{CGS-2}\label{ss:qradd_cgs2}
The CGS-2 algorithm is  simply CGS with reorthogonalization, which,
unlike CGS, is known to mitigate large cancellation errors and 
maintain numerical stability 
\cite{abdelmalek1971, daniel1976, hoffmann1989, ruhe1983}.
Reorthogonalizing the vectors in $Q$ updates the associated
orthogonal projector to include a correction matrix
$T$, although this matrix may not be explicitly formed in practice.

In Appendix 1 of \cite{LOWSYNC_GMRES}, the form of the projection
and correction matrices was derived within the context of
DCGS-2 (discussed in \cref{ss:qradd_dcgs2}), however, the matrices
remain the same and are given by
\begin{align*}
P &= I - Q_{0:m_i-2} T_{0:m_i-2,0:m_i-2} Q_{:,0:m_i-2}^T, \\
T_{0:m_i-2,0:m_i-2} &= I - L_{0:m_i-2,0:m_i-2} 
                      - L_{0:m_i-2,0:m_i-2}^T. 
\end{align*}
Here, $L$ contains the same information as in the ICWY-MGS
case, namely each row of $L$ is generated by the reorthogonalization
of the vectors within $Q$ before updating the factorization
with an additional vector. 

\cref{alg:qr_cgs2} details the process of updating the AA
$QR$ factorization with an additional vector without explicitly
forming the correction matrix. 
Reorthogonalization happens explicitly on \cref{cgs2_reorth1,cgs2_reorth2} before being added back into the
matrix $R$ on \cref{cgs2_addback}.
\begin{center}
\vspace{-2em}
\resizebox{.5\textwidth}{!}{%
\begin{minipage}{0.575\textwidth}
\begin{algorithm}[H]
    \setstretch{1.15}
    \caption{\code{QRAdd\_CGS2}}\label{alg:qr_cgs2}
    \hspace*{\algorithmicindent} \textbf{Input:} $Q$, $R$, $\Dfione$ and $m_i$ \\
    \hspace*{\algorithmicindent} \textbf{Output:} $Q$, $R$
    \begin{algorithmic}[1]
        \State $s \leftarrow Q_{:,0:m_i-2}^T * \Dfione$
        \Comment{Sync}
        \State $y \leftarrow \Dfione - Q_{:,0:m_i-2} * s$\label{cgs2_linearsum1}
        \State $z \leftarrow Q_{:,0:m_i-2}^T * y$
        \Comment{Sync}\label{cgs2_reorth1}
        \State $\Dfione \leftarrow y - Q_{:,0:m_i-2} * z$\label{cgs2_reorth2}
        \State $R_{0:m_i-2,m_i-1} \leftarrow s + z$\label{cgs2_addback}
        \State $R_{m_i-1,m_i-1} \leftarrow \| \Dfione \|_2$
        \Comment{Sync}
        \State $Q_{:,m_i-1} \leftarrow \Dfione / R_{m_i-1,m_i-1}$
    \end{algorithmic}
\end{algorithm}
\end{minipage}
}
\end{center}
Overall, this process requires three reductions to be performed
per iteration---two for the orthogonalization and reorthogonalization
of $\Dfione$ against all the vectors in $Q$ and one for the final
normalization of $\Dfione$.
The amount of computation per iteration has increased to
approximately $2 m_i n$ or $\BigO(m_i n)$. While the amount
of computation is still on 
the same order as MGS, the workload has
approximately doubled, and, fortunately, there are no additional
storage requirements.

\subsection{DCGS-2}\label{ss:qradd_dcgs2}
DCGS-2, or CGS with delayed reorthogonalization,
is based on the delayed CGS-2 algorithm introduced by 
Hernandez  et al.\ in \cite{HERNANDEZ}. As the name suggests,
the reorthogonalization of the vectors in CGS-2 is delayed
to the subsequent iteration. In his original derivation,
Hernandez notes that this process is tantamount to updating the $QR$
factorization with a CGS vector, and a deteriorated 
loss of orthogonality 
(between $\BigO(\varepsilon)$--$\BigO(\varepsilon)\kappa^2(A)$, 
that of CGS-2 and CGS) is often observed.

A stable variant of DCGS-2 was derived by
\'{S}wiryzdowicz  et al.\ \cite{LOWSYNC_GMRES} and Bielich  et al.\
\cite{DCGS2Arnoldi} by exploiting the
form of the correction matrix $T$ and introducing a normalization lag
and delayed reorthogonalization.
The symmetric correction matrix $T_{0:m_i-2,0:m_i-2}$
is the same as in the context of CGS-2, namely
\[ 
T_{0:m_i-2,0:m_i-2} = I - L_{0:m_i-2,0:m_i-2} - L_{0:m_i-2,0:m_i-2}^T.
\]
When the matrix $T_{0:m_i-2,0:m_i-2}$ is split into the pieces
$I - L_{0:m_i-2,0:m_i-2}$ and $L_{0:m_i-2,0:m_i-2}^T$, then
applied across two iterations of the DCGS-2 algorithm
coupled with lagging the normalization of the added vector, 
the resulting loss of  orthogonality is $\BigO (\varepsilon)$ in practice.
Performing DCGS-2 in this way decreases the number of reductions to 
one per iteration because the reorthogonalization is performed ``on-the-fly" 
and essentially operates a single iteration behind. The authors note 
in their paper that the final iteration of GMRES requires an additional
synchronization due to the required final normalization.

\cref{alg:qr_dcgs2} details our \code{QRAdd\_DCGS2} implementation
in which there are two synchronizations per iteration. Because
AA requires using the $QR$-factorization to solve the LSP at 
the end of \textit{every} iteration,
it is not possible to lag the normalization of the vector $\Dfione$.
\begin{center}
\vspace{-2em}
\resizebox{.5\textwidth}{!}{%
\begin{minipage}{0.575\textwidth}
\begin{algorithm}[H]
    \setstretch{1.15}
    \caption{\code{QRAdd\_DCGS2}}\label{alg:qr_dcgs2}
    \hspace*{\algorithmicindent} \textbf{Input:} $Q$, $R$, $\Dfione$ and $m_i$ \\
    \hspace*{\algorithmicindent} \textbf{Output:} $Q$, $R$
    \begin{algorithmic}[1]
        \State $R_{0:m_i-2,m_i-1} \leftarrow Q_{:,0:m_i-2}^T * \Dfione$
        \Comment{Delayed Sync}\label{dcgs2_orthodfi}
        \If{$m_i > 3$}
            \State $s \leftarrow Q_{:,0:m_i-3}^T Q_{:,m_i-2}$
            \Comment{Sync}\label{dcgs2_reortho}
            \State $Q_{:,m_i-2} \leftarrow Q_{:,m_i-2} - Q_{:,m_i-3} * s$\label{dcgs2_linearsum1}
            \State $R_{0:m_i-3,m_i-2} \leftarrow R_{0:m_i-3,m_i-2} + s$
        \EndIf
        \State $\Dfione \leftarrow \Dfione - Q_{:,0:m_i-2} * R_{0:m_i-2,m_i-1}$ \label{dcgs2_linearsum2}
        \State $R_{m_i-1,m_i-1} \leftarrow \| \Dfione \|_2$
        \Comment{Sync}
        \State $Q_{:,m_i-1} \leftarrow \Dfione / R_{m_i-1,m_i-1}$
    \end{algorithmic}
\end{algorithm}
\end{minipage}
}
\end{center}
The global reduction for the reorthogonalization of the vectors can be still 
be performed at the same time as the orthogonalization of $\Dfione$ against
$Q$. The synchronization on \cref{dcgs2_orthodfi} is lagged to
coincide with \cref{dcgs2_reortho}.
The amount of computation for DCGS-2 increases over that of CGS-2
except on the first iteration where it performs the same number 
of flops as MGS. Lagging the reorthogonalization and introducing
the operation $Q_{:,0:m_i-3}^T Q_{:,m_i-2}$, 
makes this method $\BigO(m_i n^2)$, the same as ICWY-MGS.

\section{Performance Study}\label{sec:perf_study}

In this section we detail the GPU and CPU performance of the 
\code{QRAdd} variants from \cref{sec:qradd}.
We differentiate between two performance costs:
the time required to fill the AA
iteration space, defined as ``Start-Up Iterations'' and the time
required to add an additional vector to the space once it is filled,
termed ``Recycle Iteration''.
In the case of CPU performance, a parallel strong-scaling study is performed.
With GPUs, we perform a weak-scaling study,
because a common goal is to saturate the devices, 
and thus avoid under-utilization.

Matrices and vectors containing $n$ rows are partitioned
row-wise across $p$ processes when performing
strong-scaling studies on CPUs. That is, each process (or CPU core) contains
$n/p$ contiguous rows of a given matrix or vector. For GPU
performance studies, we provide the local vector size per GPU
participating in the computation.
All tests are performed on the LLNL Lassen supercomputer \cite{lassen}, 
which has two IBM Power-9 CPUs per node each with 20 available cores,
and each CPU is connected to two Nvidia V100 GPUs.
For GPU performance tests each GPU is paired with a single MPI
rank as the host process.
Each test is performed 20 times and for 10 iterations (for a total of 200 timings); 
for each individual test the maximum time required by any single process or GPU is recorded, and the minimum time across all tests is presented.

The results presented
are based on the AA implementation within SUNDIALS\@. Parallel CPU tests are
implemented with a so-called node-local vector abstraction, called the Parallel N\_Vector. Similarly, the GPU tests are
implemented with the MPIPlusX N\_Vector where the CUDA N\_Vector
is used as the local vector for each MPI rank \cite{BALOS2021102836}.
Notably, while we use the CUDA N\_Vector as the local vector portion, this
can be switched with any N\_Vector implementation in
SUNDIALS\@. Due to this abstraction, communication between GPUs
is still performed by staging data through the host process
and not via CUDA-aware technologies. 

The low synchronization methods added to SUNDIALS leverage the fused dot product operation to perform matrix-vector products with $Q^T$. This operation enables computing the dot product of a single vector with  multiple vectors (the columns of $Q$) as a single operation requiring only one MPI call. For ICWY and DCGS-2 a new N\_Vector operation was introduced in SUNDIALS enabling delayed synchronization by separating the local reduction computation and final global reduction into separate operations. Both of these fused operations perform a one-time copy of the independent vector data pointers into a device buffer accessed by the fused operations. 
In addition to combining multiple reductions into a single call, these operations reduce the number of kernel launches in the GPU case. We discuss further in \cref{ss:startup}.

\subsection{Start-Up Iterations}\label{ss:startup}

\begin{figure*}[!ht]
  \centering
  \includegraphics{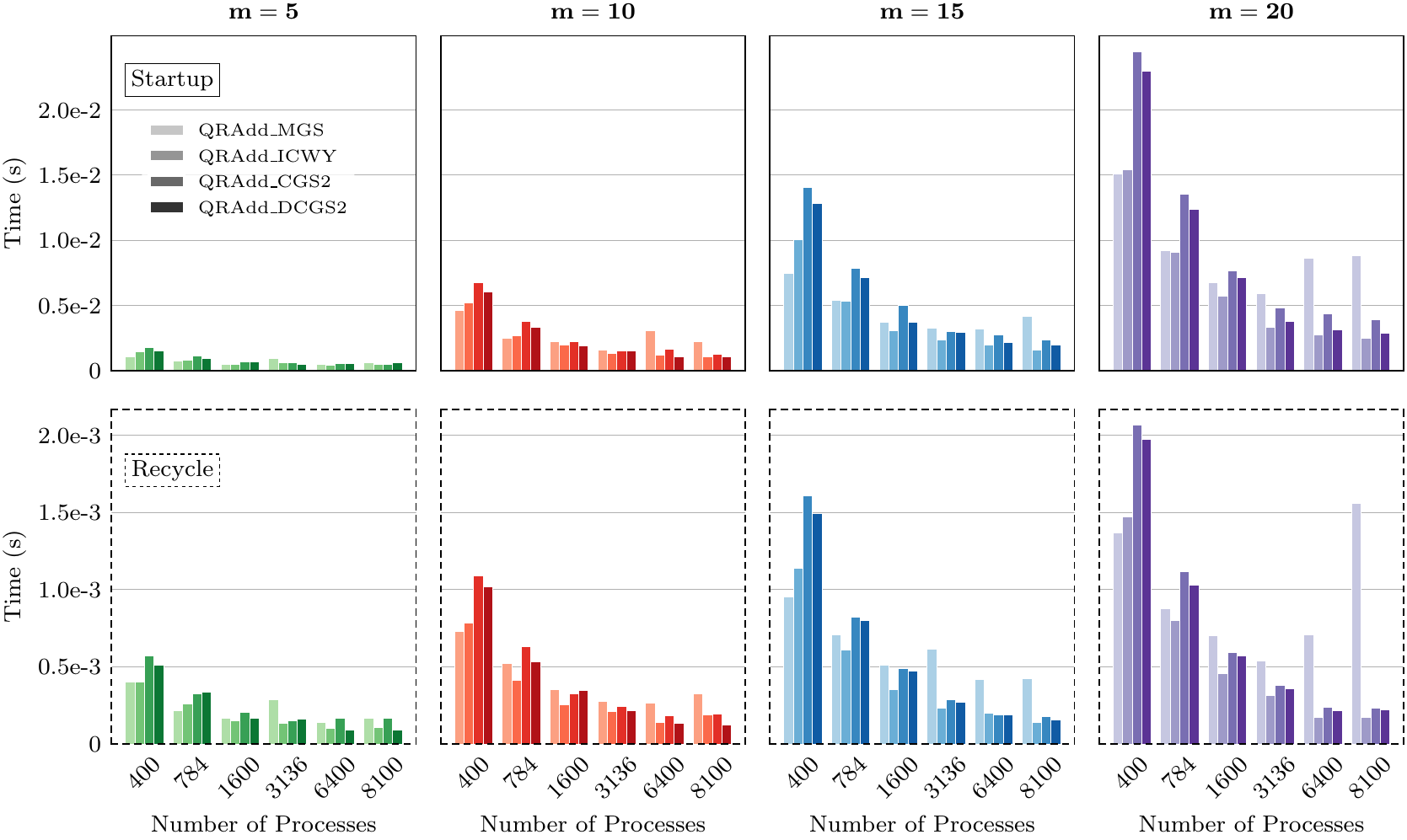}
  \caption{{\small\baselineskip=9pt
  The cumulative time for the \textbf{start-up iterations (top)} and
    \textbf{single recycle iteration (bottom)} of the four different \code{QRAdd} kernels within AA
for a vector of size $\num{1 048 576}$ on a varying number of Lassen \textbf{CPU} cores.}}\label{fig:startup_recycle_lassen_cpu}
\end{figure*}

\begin{figure*}[!ht]
  \centering
  \includegraphics{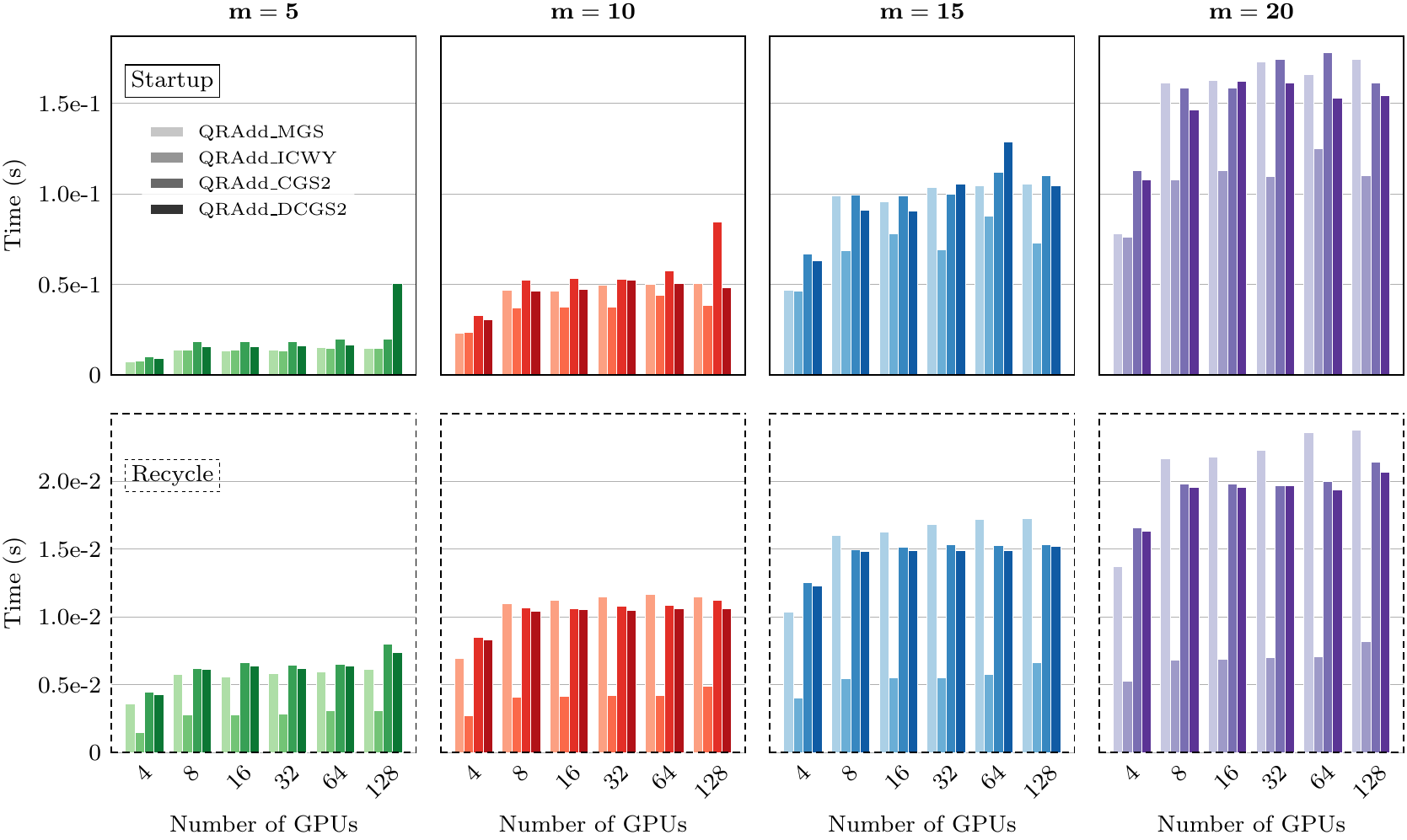}
  \caption{{\small\baselineskip=9pt
  The cumulative time for the \textbf{start-up iterations (top)} and
    \textbf{single recycle iteration (bottom)} of the four different \code{QRAdd} kernels within AA
for a local vector of size $\num{1 500 000}$ on a varying number of Lassen \textbf{GPUs}.}}\label{fig:startup_recycle_lassen_gpu}
\end{figure*}

We first consider the start-up iteration times, displayed in \cref{fig:startup_recycle_lassen_cpu,fig:startup_recycle_lassen_gpu}.
On the CPU we use a global vector of size $\num{1 048 576}$, and
for the GPU a local vector size of $\num{1 500 000}$.
The number of iterations performed is equivalent to the number of vectors, $m$,
in the AA iteration space.

The top rows of \cref{fig:startup_recycle_lassen_cpu} (CPU)
and \cref{fig:startup_recycle_lassen_gpu} (GPU) display the
time required to
fill the AA iteration space.  As part of this operation, \code{QRAdd\_MGS}
requires
$\sum_{k=1}^{m} k = \frac{m^2 + m}{2}$ total synchronizations, while
\code{QRAdd\_ICWY}, \code{QRAdd\_CGS2}, and \code{QRAdd\_DCGS}
require $2m-1$, $3m-2$, and $2m-1$, respectively.

Predicting performance based solely on the number of dot products
performed, we expect that \code{QRAdd\_ICWY} and \code{QRAdd\_DCGS2}
will become faster than \code{QRAdd\_MGS} after $m=3$ and
\code{QRAdd\_CGS2} will outperform \code{QRAdd\_MGS} after $m=6$,
for processor counts
where the synchronization imposed by dot products is the dominate
cost. This is consistent with the CPU performance seen in
\cref{fig:startup_recycle_lassen_cpu}.
The benefits of applying low synchronization
\code{QRAdd} algorithms is not observed until the processor count
reaches 1600 or higher when the global reduction costs lead to
larger bottlenecks. As $m$ increases, the disparity between the time
required for \code{QRAdd\_MGS} and the low synchronization variants
greatly increases.

Our observed GPU performance 
depends on the number of kernel calls for each
\code{QRAdd} subroutine in addition to the number of global
reductions.
Compared to \code{QRAdd\_MGS}, \code{QRAdd\_ICWY} performs 
fewer CUDA kernels each iteration for $m>3$,
while \code{QRAdd\_CGS2} and \code{QRAdd\_DCGS2} perform
four additional CUDA kernel calls for $m>1$.
Also, all three low synchronization
variants perform sequential on-CPU updates to $R$, unlike
\code{QRAdd\_MGS}.
Finally, in addition to the kernel launches and sequential updates,
 the low synchronization routines increase the amount of data
 being moved between the host and the GPU. Specifically, the 
 fused dot products employed by these routines require that at most
 $m_i-1$ values be copied to and from the GPU multiple times per 
 iteration while \code{QRAdd\_MGS} copies only a single value per individual dot product. 

All three low synchronization variants incur the cost of two kernel
launches that require data transfers to and from the host process.
However, in the case of \code{QRAdd\_ICWY}, the lower computational
requirements result in fewer overall CUDA kernel calls
than for \code{QRAdd\_CGS2} and \code{QRAdd\_DCGS2}.
Specifically, \cref{cgs2_linearsum1} in
\cref{alg:qr_cgs2} and \cref{dcgs2_linearsum1} in \cref{alg:qr_dcgs2}
correspond to additional vector updates that are not required
by \code{QRAdd\_ICWY};
hence their performance times remain higher than 
\code{QRAdd\_ICWY} independent of GPU count or $m$.
There are minimal performance gains when using \code{QRAdd\_CGS2}
or \code{QRAdd\_DCGS2} for filling the AA iteration space, as
\code{QRAdd\_DCGS2} is only slightly faster than \code{QRAdd\_MGS}
when more than 8 GPUs are used and only for $m=15$ and larger.

\subsection{Recycle Iterations}\label{ss:recycle}

We next consider the cost of the \code{QRAdd} kernel once the AA iteration space is
filled, namely the time required to orthogonalize a single vector against $m-1$
vectors.
Assuming a large iteration count is required for convergence, this cost 
reflects the expected run-time of the \code{QRAdd} kernel for the majority of a given solve.
We label these times as the ``Recycle Iteration'' time in
\cref{fig:startup_recycle_lassen_cpu,fig:startup_recycle_lassen_gpu}.
The CPU and GPU timings presented are for the same global
and local vector size as presented in \cref{ss:startup},
and the results are for a \emph{single} iteration.
Note that this time does not include the additional
synchronization introduced by ICWY into \code{QRDelete}, first mentioned
in \cref{ss:qradd_icwy} and discussed further in \cref{sec:experiments}.

Once the AA iteration space is filled,
the per iteration cost is greatly decreased
by using one of the low synchronization orthogonalization algorithms
when operating on CPUs (see the bottom row of \cref{fig:startup_recycle_lassen_cpu}).
While the ``Startup Iterations'' did not result in performance
gains on CPUs until after 1600 processes, for $m > 5$, we observe
performance gains with \code{QRAdd\_ICWY} with process
counts as low as 784.
\code{QRAdd\_MGS} is still faster at smaller scales for all values of
$m$ due to the reduced synchronization costs of performing an
\code{MPIAll\_Reduce} on a small, closed-set number of processes.
The performance gains of the low synchronization algorithms at larger
scales ranges from 2--8$\times$ speedup for $m=10$, $15$, and $20$ at
$\num{8 100}$ processes.
With this drastic speedup for each iteration, we 
expect much larger performance gains for test
problems that run to convergence.

Unlike the CPU performance for AA, \code{QRAdd\_ICWY} demonstrates
performance gains independent of the number of GPUs participating
in the computation or the size of the AA iteration space, $m$.
This is seen in the bottom row of
\cref{fig:startup_recycle_lassen_gpu}, which displays the performance
for a single ``Recycle Iteration'' on multiple GPUs.
While \code{QRAdd\_CGS2} and \code{QRAdd\_DCGS2} do not see the
same performance improvements as \code{QRAdd\_ICWY}, they
do begin exhibiting faster performance than
\code{QRAdd\_MGS} beginning at 8 GPUs and for $m=10$
and larger, though gains are modest.

\section{Numerical Experiments}\label{sec:experiments}

In this section, we highlight
the strong-scaling parallel efficiency of standard AA compared with
AA with low synchronization orthogonalization for test problems run to
convergence. The example
problems provided are not exhaustive, but the selected
tests do stress the performance of low synchronization AA\@. In addition,
application specific tests are used to
demonstrate the benefits of low synchronization
AA for both distributed CPU and GPU computing environments.
All experiments are performed on the LLNL Lassen supercomputer \cite{lassen}, with the same setup as described in \cref{sec:perf_study}.
To account for machine variability, each run is executed
10 times and we report the minimum.

\subsection{Anisotropic 2D Heat Equation + Nonlinear Term}
This test problem highlights the performance of all
AA variants for various $m$ and iterations to convergence.
We consider a steady-state 2D heat equation
with an additional nonlinear term $c(u)$,
\begin{align*}
  u_{xx} + u_{yy} + c(u) &= f \quad \text{in} \quad \mathcal{D} = [0,1] \times [0,1] \\
  u &= 0 \quad \text{on} \quad \partial \mathcal{D}.
\end{align*}
The chosen analytical solution is
\begin{equation*}
    u_{\text{exact}} = u(x,y) = \sin^2(\pi x) \sin^2(\pi y),
\end{equation*}
hence, the static term $f$ is defined as follows
\begin{equation*}
\begin{aligned}
    f(x,y) = \,&2\pi^2(\cos^2(\pi x) - \sin^2(\pi x)) \sin^2 (\pi y) \\
        &+ 2\pi^2(\cos^2(\pi y) - \sin^2(\pi y)) \sin^2 (\pi x) \\
        &+ c(u_{\text{exact}}).
\end{aligned}
\end{equation*}
The spatial derivatives are computed using second-order
centered differences, with the data distributed over
$\num{1024} \times \num{1024}$ points on a uniform spatial grid,
resulting in a system of equations of size
$\num{1048576} \times \num{1048576}$.
The Laplacian term is implemented as a matrix-vector product giving the the algebraic system as
\begin{equation*}\label{eq:heat2d_mat_vec}
  A \vec{u} + c(\vec{u}) = \vec{b}.
\end{equation*}
where $\vec{u}$ denotes the discrete vector of unknowns. Solving for $\vec{u}$ results in the following FP formulation
\begin{equation*}
  \vec{u} = G(\vec{u}) = A^{-1} (\vec{b} - c(\vec{u})).
\end{equation*}
We use
the SUNDIALS PCG solver to solve the
linear system with the hypre PFMG preconditioner
performing two relaxation sweeps per iteration.
Both the FP nonlinear solver and the PCG linear solver
set a stopping criteria with a tolerance
of $10^{-10}$. A zero vector is used as
the starting guess in all cases.

\subsubsection{Nonlinear Term 1}\label{ss:heat2d_term1}

\begin{figure}[!ht]
  \includegraphics{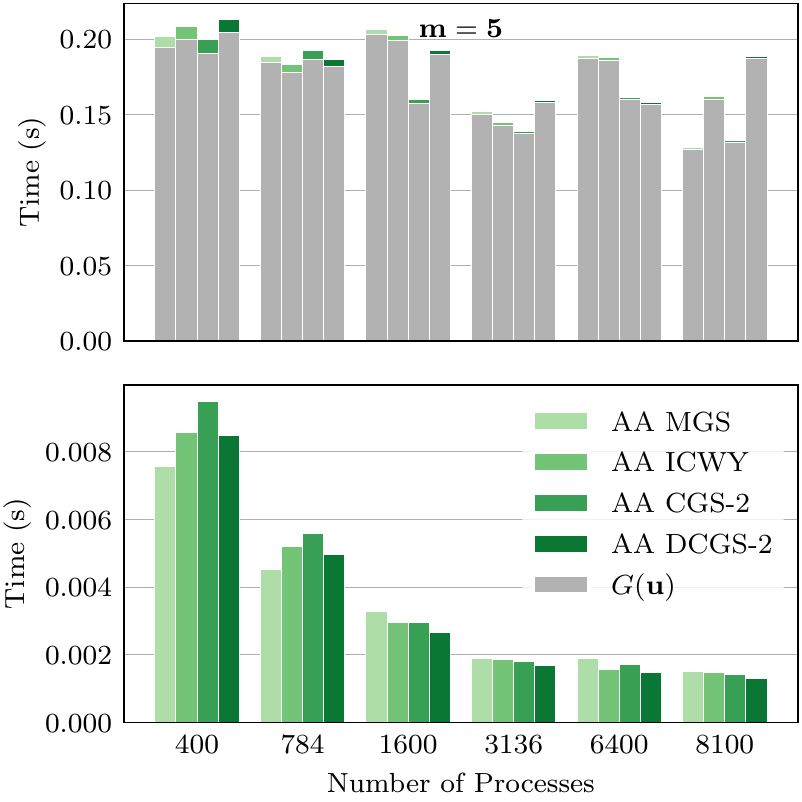}
      \caption{{\small\baselineskip=9pt
        Time to convergence for \textbf{Heat 2D + Nonlinear Term 1} (11 iterations for all cases),
        using FP + AA with $m=5$, including the function evaluation, $G(\vec{u})$, (top) and only time in AA\@ (bottom).
       }}\label{fig:heat2d_c16}
\end{figure}

As a first example, consider the nonlinear reaction term
\begin{equation*}
    c(u) = u + ue^u + ue^{-u} + (u - e^u)^2.
\end{equation*}
In this case, AA exhibits rapid
convergence when $m=5$, requiring
11 iterations to converge for all variants.
\Cref{fig:heat2d_c16} displays the overall time to convergence, split into $G(\vec{u})$ evaluation time
and time spent in AA\@. In general,
the $G(\vec{u})$ performance is volatile, most likely
due to the sparse matrix operations required for the linear solve;
thus we focus on the AA performance
in the bottom of \cref{fig:heat2d_c16}.

The time spent in AA meets expectations, based on the results
in \cref{sec:perf_study}.
For such a small $m$, we observe minimal performance improvements over AA
with MGS, particularly at small process counts.
However, the observed improvements begin around 1600 processes.

\subsubsection{Nonlinear Term 2}

As a second example, consider
\begin{equation*}
    c(u) = 100 \cdot (u - u^2).
\end{equation*}
With $m=10$,
AA requires more time to converge than for the $c(u)$ in
\cref{ss:heat2d_term1}, better highlighting potential
performance benefits of the low synchronization orthogonalization algorithms.

\begin{figure}[!ht]
  \includegraphics{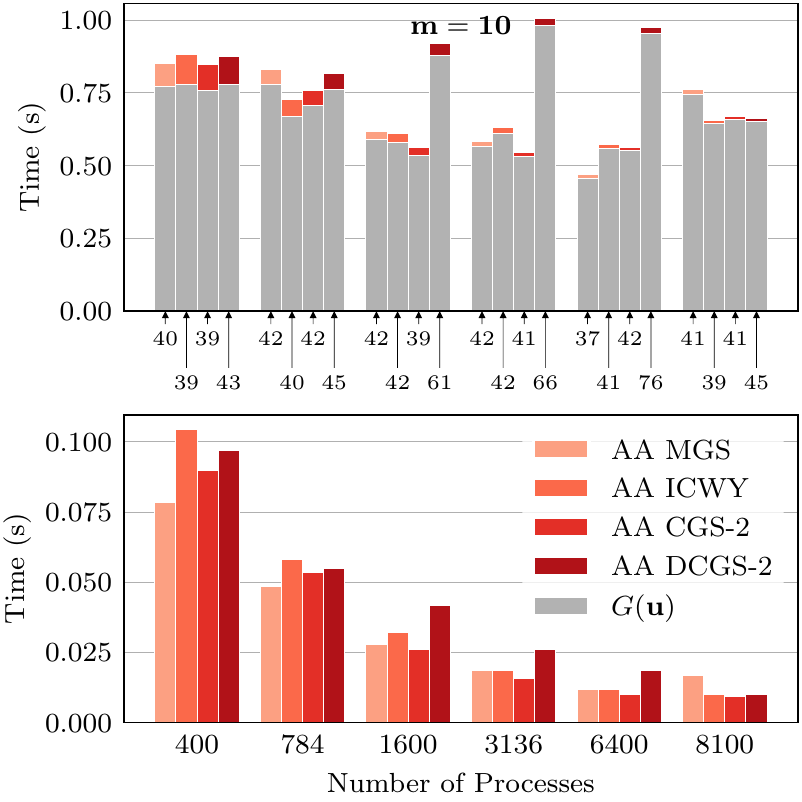}
      \caption{{\small\baselineskip=9pt
        Time to convergence for \textbf{Heat 2D + Nonlinear Term 2},
        using FP + AA with $m=10$, including the function evaluation, $G(\vec{u})$, and the number
        of iterations to convergence (top) and only the time spent in AA (bottom).}}\label{fig:heat2d_c3_perf}
\end{figure}

\Cref{fig:heat2d_c3_perf} displays the overall timing results,
with the
number of iterations required to
converge listed beneath each bar.
In most cases, the number of iterations
is between 27 and 42, depending on the orthogonalization
method and processor count. The three exceptions are
for AA with DCGS-2, which requires more iterations to
converge, particularly for 1600--6400 processes.
This is expected
since our \code{QRAdd\_DCGS2} cannot lag the normalization
and exhibits the same instability as Hernandez's
original algorithm.

In most cases, AA with DCGS-2 results in degraded
performance in comparison to the other variants due to the high number of
iterations required to converge. AA with CGS-2 performs best,
starting at 1600 processes and continues to be the fastest
through $\num{8100}$ processes. Although \code{QRAdd\_ICWY} only
requires two global synchronizations per iteration, there is
an additional synchronization in \code{QRDelete} to update
the matrix $T$, resulting in an implementation that is more comparable to
AA with CGS-2 at higher process counts.

\subsection{2D Bratu Problem}
We next consider an example from \cite{WALKER_NI} (and
originating from \cite{BRATU}).
The Bratu problem is a nonlinear PDE boundary value problem
defined as:
\begin{align*}
  u_{xx} + u_{yy} + \lambda e^u &= 0 \quad \text{in} \quad \mathcal{D} = [0,1] \times [0,1] \\
  u & = 0 \quad \text{on} \quad \partial \mathcal{D}.
\end{align*}
We again solve the problem with centered differencing,
resulting in a system of the form
\begin{equation*}
  A \vec{u} + \lambda e^{\vec{u}} = 0.
\end{equation*}
The associated FP function is then
\begin{equation*}
  G(\vec{u}) = A^{-1} \left(- \lambda e^{\vec{u}}\right) = \vec{u}.
\end{equation*}
We use a uniform grid of
$1024 \times 1024$ points, which
results in a system of equations of size
$\num{1048576} \times \num{1048576}$. We select $\lambda=6.7$ for testing,
as it is close to the theoretical critical point,
as discussed in \cite{BRATU_IC}, and is a difficult 
problem to solve.
The SUNDIALS PCG solver is applied to perform the
linear system solve with the hypre PFMG preconditioner
performing two relaxation sweeps. Both the nonlinear
and linear solver employed a tolerance of $10^{-10}$.
A zero vector is set as the starting guess.

\begin{figure}[!ht]
  \includegraphics{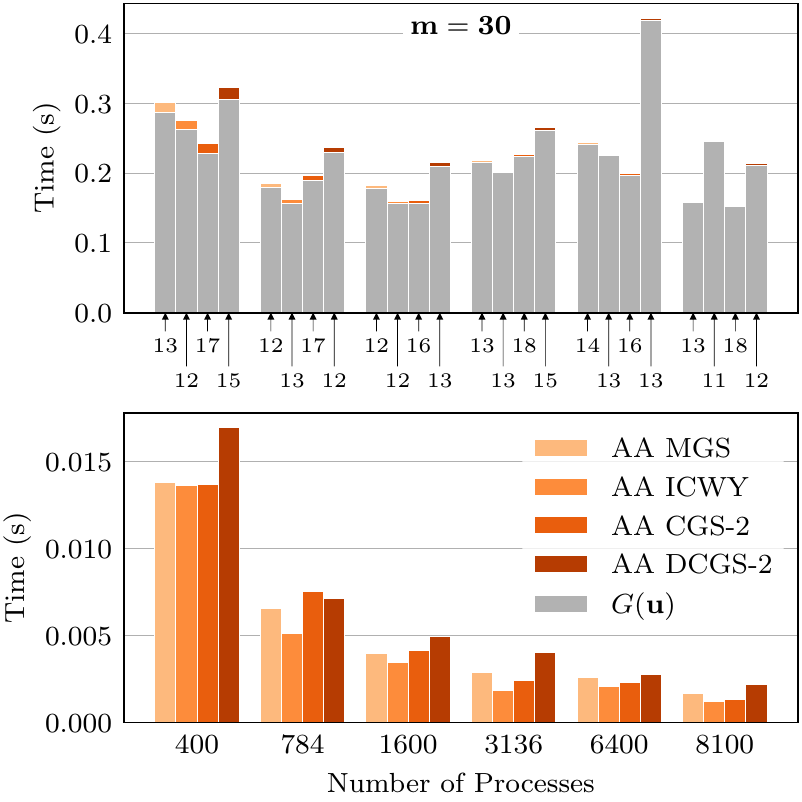}
      \caption{{\small\baselineskip=9pt
      Time to convergence for \textbf{2D Bratu},
        using FP + AA with $m=30$, including the function evaluation, $G(\vec{u})$, and the number
        of iterations to convergence (top) and only the time spent in AA (bottom).}}\label{fig:bratu2d_perf}
\end{figure}

For this example problem, $m = 30$, and
the number of iterations required to converge for all variants of AA is
less than 30; hence the timing results reflect the performance benefits of
the \code{QRAdd} subroutines for the startup iterations.
The strong-scaling timing and convergence results are presented in
\cref{fig:bratu2d_perf}.
For this test case, we observe improvements with ICWY from the very
beginning with 400 processes (although only slightly in this case). AA
with ICWY
continues to outperform up to 8100 processes.
This is consistent with the results in \cref{ss:startup} in which
\code{QRAdd\_ICWY} performed best for $m=20$ as it only requires
two dot products per iteration, one fewer than \code{QRAdd\_CGS2}
and has approximately the same amount of computation as \code{QRAdd\_DCGS2}.

\subsection{Expectation-Maximization Algorithm for Mixture Densities}

For this example, we consider a variation of the expectation-maximization
test problem presented in \cite{WALKER_NI}.
Consider a mixture density of three univariate
normal densities with a mixture density given by
$p(x) = \sum_{i=1}^3 \alpha_i p_i(x|\mu_i,\sigma_i)$, with
\begin{equation*}
    p_i(x | \mu_i,\sigma_i) = \frac{1}{\sqrt{2\pi} \sigma_i}
        e^{-(x-\mu_i)^2 / (2\sigma_i^2)}, \qquad 1 \leq i \leq 3
\end{equation*}
Mixture proportions $\{\alpha_i\}_{i=1}^3$ are non-negative and sum to one.
The mixture proportions and variances are assumed to be  known and
the means $\{\mu_i\}_{i=1}^3$ are estimated from a set of unlabeled samples $\{x_k\}_{k=1}^N$,
or samples of unknown origin.
Determining the unknown mean distribution parameters is given by
the FP function
\begin{equation*}
    G(\mu_i) = \mu_i =
    \frac{\sum_{k=1}^N x_k \frac{\alpha_i p_i(x_k | \mu_i, \sigma_i)}{p(x_k)} }{
          \sum_{k=1}^N \frac{\alpha_i p_i(x_k | \mu_i, \sigma_i)}{p(x_k)} },
          \qquad 1 \leq i \leq 3
\end{equation*}
with current mean estimations $\{\mu_i\}_{i=1}^3$ being applied alongside the known
mixture proportions and variances to determine the subsequent estimations
until convergence.
We keep the same mixture proportions and variances as the original test case,
$(\alpha_1, \alpha_2, \alpha_3) = (0.3, 0.3, 0.4)$ and
$(\sigma_1, \sigma_2, \sigma_3) = (1, 1, 1)$.
We generated \num{100 000} samples for the mean distribution set
$(\mu_1=0, \mu_2=0.5, \mu_3=1.0)$,
corresponding to a poorly separated mixture, and used the same 
AA parameter of $m=3$ as Walker and Ni \cite{WALKER_NI}. 

For our test case, we estimate a single set of mean distribution
parameters redundantly for every entry in a global vector, where 
the vector takes the form
\begin{equation*}
    \vec{u} = \begin{bmatrix} \{\mu_i\}_{i=1}^3 & \cdots & \{\mu_i\}_{i=1}^3 \end{bmatrix}.
\end{equation*}
We do this to simulate a function that requires no communication
other than that imposed by AA.
The resulting FP function to be solved is then given by
\begin{equation*}
    G(\vec{u}) = \begin{bmatrix} \{G(\mu_i)\}_{i=1}^3 & \cdots & \{G(\mu_i)\}_{i=1}^3 \end{bmatrix}.
\end{equation*}
Because $m$ is small for this test,
we expect to see only modest improvements in performance of the low
synchronization routines over MGS. For $m=3$, ICWY and DCGS-2 reduce
the number of synchronizations per iteration by one over MGS, with
ICWY gaining an additional synchronization after the space is filled.

The test is performed as a weak-scaling study with each GPU operating on a
local vector size of \num{1 500 000} values.
Tests were run with a tolerance of $10^{-8}$, and each AA version requires 21 iterations to converge, independent of the
orthogonalization subroutine used.
\begin{figure}[!ht]
  \includegraphics{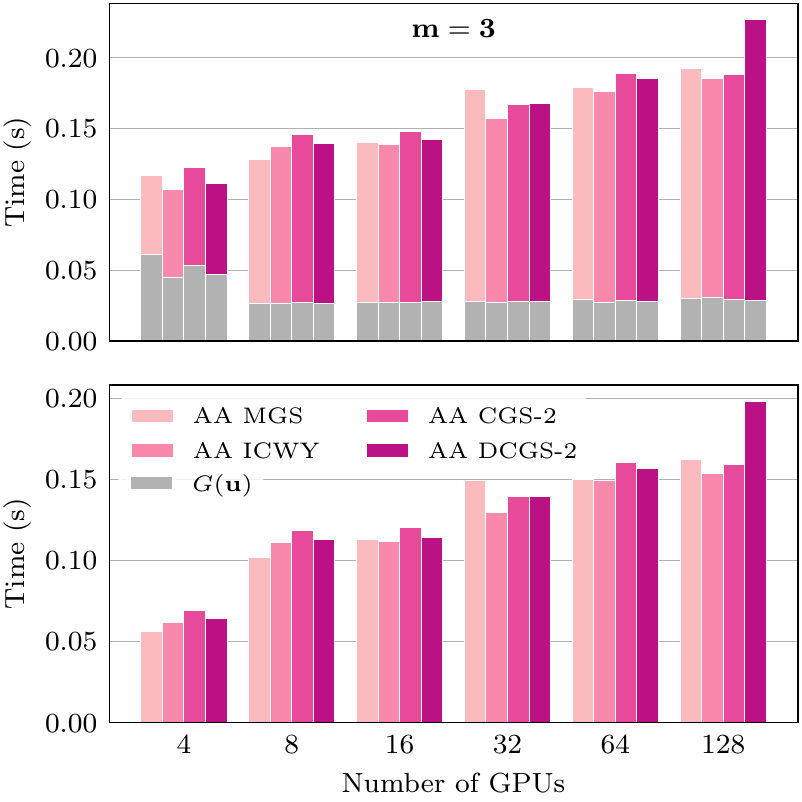}
      \caption{{\small\baselineskip=9pt
      Time to convergence for \textbf{Expectation-Maximization} (21 iterations for all cases),
        using FP + AA with $m=3$
        and a local vector size of \num{1 500 000} for each GPU, including 
        the function evaluation, $G(\vec{u})$, (top) and only the time spent in AA (bottom). }}\label{fig:em_perf}
\end{figure}
The results of these tests are presented in  \cref{fig:em_perf} and are consistent with the
results of the weak-scaling study since no communication is required for $G(\vec{u})$. The function
evaluation performance remains the same, independent of GPU count or AA variant
(differing only slightly for runs on a single node with 4 GPUs).
In addition, the time spent in $G(\vec{u})$ is lower for this example since $G(\vec{u})$ performance depends on the number of samples in the
mixture, while the performance of AA scales with the number of values in the 
global vector.

Results are consistent with observations in the GPU weak-scaling study in
\cref{sec:perf_study} where all of the \code{QRAdd} subroutines performed similar for
$m=5$. Additionally, because $m$ is small, we do not incur
a large overhead for data movement and multiple kernel launches in comparison to those observed
with larger values of $m$, as presented in \cref{ss:recycle}. Because the majority of the
computation has moved to the GPU, we observe the expected cross-over point for 
ICWY with its performance being slightly faster than that of
MGS for GPU counts of 16 and larger. There are no consistent improvements
observed for CGS-2 or DCGS-2.  

\section{Conclusions}\label{sec:conclusions}

Anderson Acceleration (AA) is an efficient method for accelerating
the convergence of fixed point solvers, but faces performance
challenges in parallel distributed computing environments
mainly due to the number of global synchronizations
per iteration which is dependent upon the size of the AA iteration space.
In this paper, we introduced low synchronization orthogonalization
subroutines into AA which effectively
reduce the number of global synchronizations to a constant number
per iteration independent of the size of the AA iteration
space. We presented a performance study that demonstrated
the improved strong-scalability of AA with these low synchronization
\code{QRAdd} subroutines when performed in a CPU-only parallel
environment, as well as demonstrated performance and implementation
concerns for these subroutines when operating in a multi-GPU
computing environment.
Furthermore, our numerical results display a realistic
picture of the expected performance of AA
in practice that matches the predictions of our performance
analysis in \cref{sec:perf_study} and suggests the use of ICWY for
large values of $m$ when operating in a CPU-only parallel environment and 
as the default method for distributed GPU computing.
Overall, this paper provides a comprehensive study of low synchronization
orthogonalization routines within AA and their parallel performance benefits.

The software used to generate the results in this paper
will be available in a future release of SUNDIALS.

\section*{Acknowledgments}

Support for this work was provided through the Scientific Discovery through Advanced Computing (SciDAC) project ``Frameworks, Algorithms and Scalable Technologies for Mathematics (FASTMath),'' funded by the U.S. Department of Energy Office of Advanced Scientific Computing Research.

This work was performed under the auspices of the U.S. Department of Energy by Lawrence Livermore National Laboratory under contract DE-AC52-07NA27344. Lawrence Livermore National Security, LLC. LLNL-PROC-827159.

This material is based in part upon work supported by the Department of Energy,
National Nuclear Security Administration, under Award Number \textit{DE-NA0003963}.

\bibliographystyle{siamplain}
\bibliography{paper-lowsyncanderson.bib}

\end{document}